\magnification=1200
\hsize=12.5cm
\vsize=19.5cm
\parskip 3pt plus 1pt minus 1pt
\parindent=6.4mm

\let\sl\it

\font\eightrm=cmr8
\font\eightbf=cmbx8
\font\eightsy=cmsy8
\font\eighti=cmmi8
\font\eightit=cmti8
\font\eighttt=cmtt8
\font\eightsl=cmsl8

\font\sixrm=cmr6
\font\sixbf=cmbx6
\font\sixsy=cmsy6
\font\sixi=cmmi6

\font\fourteenbf=cmbx10 at 14.4pt
\font\twelvebf=cmbx10 at 12pt

\font\tenmsa=msam10
\font\sevenmsa=msam7
\font\fivemsa=msam5
\newfam\msafam
  \textfont\msafam=\tenmsa 
  \scriptfont\msafam=\sevenmsa
  \scriptscriptfont\msafam=\fivemsa

\font\tenmsb=msbm10
\font\eightmsb=msbm8
\font\sevenmsb=msbm7
\font\fivemsb=msbm5
\newfam\msbfam
  \textfont\msbfam=\tenmsb
  \scriptfont\msbfam=\sevenmsb
  \scriptscriptfont\msbfam=\fivemsb
\def\Bbb{\fam\msbfam\tenmsb}

\catcode`\@=11
  \def\eightpoint{%
  \textfont0=\eightrm \scriptfont0=\sixrm \scriptscriptfont0=\fiverm
  \def\rm{\fam\z@\eightrm}%
  \textfont1=\eighti \scriptfont1=\sixi \scriptscriptfont1=\fivei
  \def\oldstyle{\fam\@ne\eighti}%
  \textfont2=\eightsy \scriptfont2=\sixsy \scriptscriptfont2=\fivesy
  \textfont\itfam=\eightit
  \def\it{\fam\itfam\eightit}%
  \textfont\slfam=\eightsl
  \def\sl{\fam\slfam\eightsl}%
  \textfont\bffam=\eightbf \scriptfont\bffam=\sixbf
  \scriptscriptfont\bffam=\fivebf
  \def\bf{\fam\bffam\eightbf}%
  \textfont\ttfam=\eighttt
  \def\tt{\fam\ttfam\eighttt}%
  \def\Bbb{\fam\msbfam\eightmsb}%
  \textfont\msbfam=\eightmsb
  \def\Cal{\fam\Calfam\eightCal}%
  \textfont\Calfam=\eightCal
  \abovedisplayskip=9pt plus 2pt minus 6pt
  \abovedisplayshortskip=0pt plus 2pt
  \belowdisplayskip=9pt plus 2pt minus 6pt
  \belowdisplayshortskip=5pt plus 2pt minus 3pt
  \smallskipamount=2pt plus 1pt minus 1pt
  \medskipamount=4pt plus 2pt minus 1pt
  \bigskipamount=9pt plus 3pt minus 3pt
  \normalbaselineskip=9pt
  \setbox\strutbox=\hbox{\vrule height7pt depth2pt width0pt}%
  \let\bigf@ntpc=\eightrm \let\smallf@ntpc=\sixrm
  \normalbaselines\rm}
\catcode`\@=12

\def\bC{{\Bbb C}}

\def\bP{{\Bbb P}}
\def\bQ{{\Bbb Q}}
\def\bR{{\Bbb R}}
\def\bZ{{\Bbb Z}}

\def\bB{{\Bbb B}}

\font\tenCal=eusm10
\font\eightCal=eusm8
\font\sevenCal=eusm7
\font\fiveCal=eusm5
\newfam\Calfam
  \textfont\Calfam=\tenCal
  \scriptfont\Calfam=\sevenCal
  \scriptscriptfont\Calfam=\fiveCal
\def\Cal{\fam\Calfam\tenCal}

\def\cC{{\Cal C}}

\def\det{{\mathop{\rm det}\nolimits}}

\def\dim{\mathop{\rm dim}\nolimits}

\def\Vol{\mathop{\rm Vol}\nolimits}

\def\ddbar{\partial\overline\partial}
\def\\{\hfil\break}
\let\wt\widetilde
\let\wh\widehat

\let\ep\varepsilon

\def\hexnbr#1{\ifnum#1<10 \number#1\else
 \ifnum#1=10 A\else\ifnum#1=11 B\else\ifnum#1=12 C\else
 \ifnum#1=13 D\else\ifnum#1=14 E\else\ifnum#1=15 F\fi\fi\fi\fi\fi\fi\fi}
\def\msatype{\hexnbr\msafam}
\def\msbtype{\hexnbr\msbfam}
\mathchardef\compact="3\msatype62
\mathchardef\smallsetminus="2\msbtype72   
\mathchardef\subsetneq="3\msbtype28
\def\buildo#1\over#2{\mathrel{\mathop{\null#2}\limits^{#1}}}
\def\buildu#1\under#2{\mathrel{\mathop{\null#2}\limits_{#1}}}

\def\square{{\hfill \hbox{
\vrule height 1.453ex  width 0.093ex  depth 0ex
\vrule height 1.5ex  width 1.3ex  depth -1.407ex\kern-0.1ex
\vrule height 1.453ex  width 0.093ex  depth 0ex\kern-1.35ex
\vrule height 0.093ex  width 1.3ex  depth 0ex}}}

\let\Item\item
\def\item#1{\Item{$\rlap{\hbox{#1}}\kern\parindent\kern-5pt$}}

\def\today{\ifcase\month\or
January\or February\or March\or April\or May\or June\or July\or August\or
September\or October\or November\or December\fi \space\number\day,
\number\year}

\null
\vskip 3cm
\centerline{\fourteenbf Une g\'en\'eralisation du th\'eor\`eme de }

\bigskip
\centerline{\fourteenbf Kobayashi-Ochiai}

\null
\vskip 20pt
\centerline{\rm  Fr\'ed\'eric Campana et Mihai P\u aun}

\null

\medskip

\vskip 30pt

\noindent{\twelvebf 0. Introduction} 

\noindent Soit $X$ une vari\'et\'e k\"ahl\'erienne
compacte $X$, de dimension $n$, munie d'une m\'etrique $\omega$, et  une application non d\'eg\'en\'er\'ee $\varphi:\bC^n\mapsto X$.

Dans cet article, nous  \'etablissons des relations entre la
croissance de $\varphi$ (mesu-r\'ee par le degr\'e moyen ou la fonction caract\'eristique) et la positivit\'e du fibr\'e canonique $K_X$ (mesur\'ee par sa pseudo-effectivit\'e, et sa dimension num\'erique). Le principe g\'en\'eral \'etant que la croissance de $\varphi$ augmente avec la positivit\'e de $K_X$.

\

\noindent Afin d'introduire notre premier r\'esultat, on rappelle la notion de 
degr\'e moyen (voir par exemple le travail de K. Kodaira [9]) si $\varphi: \Delta_r\rightarrow X$
est une application holomorphe, son d\'egr\'e moyen est par d\'efinition
$$\deg(\varphi\vert \Delta_r):= \int_{\Delta_r}\varphi^*\omega^n/\int_X\omega^n.$$
Dans ce contexte, on a le th\'eor\`eme suivant (essentiellement optimal, par l'exemple des tores complexes):

\medskip
\noindent {\bf Th\'eor\`eme 1.} {\sl Soit $(X, \omega)$ une vari\'et\'e complexe
compacte de dimension $n$. S'il existe une application holomorphe 
$\varphi:\bC^n\rightarrow X$ non d\'eg\'en\'er\'ee,
et telle que

$$\lim_{r\mapsto \infty}1/r^{2n}\deg(\varphi\vert \Delta_r)= 0,\leqno (1) $$

\noindent alors le fibr\'e canonique de $X$ n'est pas pseudo-effectif.}

\medskip

\noindent Il convient de rappeller ici que sous les m\^emes hypoth\`eses,
Kodaira montre dans [9] que le tous les plurigenres de $X$
sont nuls;  le th\'eor\`eme 1 peut \^etre vu comme une am\'elioration
de son r\'esultat. 

\noindent En particulier, si la vari\'et\'e $X$ est projective, cette am\'elioration, combin\'ee avec le th\'eor\`eme principal de  [3] r\'esolvent la conjecture d'unir\'eglage dans cette situation:

\medskip

\noindent {\bf Corollaire.} {\sl Soit $X$ une vari\'et\'e projective $n$-dimensionelle, et 
$\varphi:\bC^n\to X$ une application holomorphe non  
d\'eg\'en\'er\'ee, et qui v\'erifie la condition (1). Alors $X$ est unir\'egl\'ee.}

\noindent On remarque que, d\`es que $\dim(X)\geq 2$, il existe des exemples de vari\'et\'es $X$ et d'applications holomorphes $\varphi:\bC^n\to X$
non-d\'eg\'en\'er\'ees de volume fini ((Fatou-Bieberbach) dont l'image a un compl\'ementaire contenant un ouvert non-vide.

La d\'emonstration du th\'eor\`eme 1 ne fournit aucun rapport entre
les courbes rationnelles (produites par la th\'eorie de Mori) et les images des droites
complexes par l'application $\varphi$. Intuitivement, on s'attend \`a ce que ``beaucoup'' parmi les images de droites  par $\varphi$ se compactifient, mais ceci est faux en g\'en\'eral: voir un contre-exemple \`a la fin du chapitre 1, dans lequel aucune courbe alg\'ebrique de $\bC^n$ n'a une image alg\'ebrique par $\varphi$).  Nous remercions TC Dinh et N. Sibony pour nous avoir indiqu\'e la r\'ef\' erence [12], dans laquelle on trouve une r\'eponse positive \`a la question sur les domaines de Fatou-Bieberbach que nous nous posions.

\vskip 15pt
Dans la th\'eorie de Nevanlinna \'equidimensionelle, developp\'ee entre autres 
par Griffiths, King, Stoll..., on d\'efinit la {\sl fonction caract\'eristique} $T_\omega(\varphi, .)$
d'une application holomorphe (substitut du degr\'e dans le cadre compact) par l'\'egalit\'e:

$$T_\omega(\varphi, r):= \int_0^r{{dt}\over{t^{2n-1}}}\int_{\vert z \vert< t}\varphi^*\omega
\wedge \omega_{euc}^{n-1},$$
o\`u $\omega $ est une m\'etrique hermitienne sur $X$, et 
$\omega_{euc}:= i\ddbar\vert z\vert^2$
est la m\'etrique canonique sur $\bC^n$ (voir par exemple [7] pour une interpr\'etation de cette
quantit\'e, et certaines de ses propri\'et\'es).

Consid\'erons les donn\'ees suivantes: $X$ est une vari\'et\'e k\"ahl\'erienne compacte de fibr\'e canonique $K_X$ est pseudo-effectif, et $\varphi:\bC^n\to X$ est une 
application holomorphe non-d\'eg\'en\'er\'ee. Le r\'esultat suivant 
montre que la fonction caract\'eristique 
de $\varphi$ croit d'autant plus vite que 
la dimension num\'erique de $K_X$ est grande.

\medskip

\noindent {\bf Th\'eor\`eme 2.} {\sl Soit $(X, \omega)$ une vari\'et\'e k\"ahl\'erienne
compacte, de dimension $n$, de fibr\'e canonique $K_X$ pseudo-effectif.
Consid\'erons une application holomorphe $\varphi:\bC^n\mapsto X$ non d\'eg\'en\'er\'ee.
Alors il existe une constante positive $C> 0$ telle que 

$$\displaystyle T(\varphi, r)^{1-\nu/n}\geq Cr^{2},$$ 
pour tout $r> 0$, o\`u $\nu= \nu(K_X)$
est la dimension num\'erique de $K_X$.
 
}
\medskip

\noindent Le th\'eor\`eme pr\'ec\'edent 
g\'en\'eralise ainsi le r\'esultat classique suivant, d\^u \`a Kodaira, Griffiths, 
Kobayashi--Ochiai.
\medskip
\noindent {\bf Th\'eor\`eme {\rm ([9], [7], [8])}} {\sl Soit $X$ une vari\'et\'e de dimension $n$,
de type g\'e\-n\'e\-ral. Alors pour toute application $\varphi:\bC^n\mapsto X$, on a: 
$\varphi^*(\omega^n)= 0$, et $\varphi$ est donc d\'eg\'en\'er\'ee.}

\medskip Si le fibr\'e canonique $K_X$ est seulement situ\'e sur le bord du c\^one 
{\it big}, il semble beaucoup plus d\'elicat d'analyser les cons\'equences
du fait que la pseudo-forme volume de Kobayashi de $X$ soit d\'eg\'en\'er\'ee
(\'evidemment, on ne peut pas esp\'erer un \'enonc\'e analogue au th\'eor\`eme 
pr\'ecedent, comme le montre le cas des tores). 
N\'eanmoins, on voudrait proposer le probl\`eme suivant.
\medskip
\noindent {\bf Conjecture 1.} {\sl 
Soit $(X, \omega)$ une vari\'et\'e k\"ahl\'erienne compacte
telle que $K_X$ est pseudo-effectif. 
Supposons que $X$ est rev\^etue par $\bC^n$.

Alors la dimension num\'erique de $K_X$ est zero.(Lorsque $X$ est projective, la condition de pseudo-effectivit\'e est superflue, puisque $K_X$ est nef).}

\medskip 

Cette conjecture g\'en\'eralise
Kobayashi-Ochiai et implique une c\'el\`ebre conjecture d'Iitaka, par des arguments
standard.

Pour les courbes enti\`eres $\varphi: \bC\mapsto X$, on dispose du lemme de
``re\-pa\-ram\'e\-trisation'' de Brody: \'etant donn\'ee $\varphi$, non-constante, il existe 
une application $\wh \varphi :\bC\mapsto X$, telle que 
$\sup_{\bC}\vert \wh \varphi^{\prime}(t) \vert = 1$. Le proc\'ed\'e de Brody ne donne que 
des r\'esultats partiels en plusieurs variables: si $\varphi: \bC^m\mapsto X$ est
une application holomorphe, alors il existe une suite $\varphi_k: \Delta_k\mapsto X$
(ici, $\Delta_k$ d\'esigne le polydisque de rayon \'egal \`a $k$) telle que
$\displaystyle \varphi_k^*\omega^m\leq {{d\lambda}\over {\prod_j (1-\vert t_j/k\vert ^2)^2}}$,
avec \'egalit\'e \`a l'origine. On peut montrer que si, de plus, les images des
applications $\varphi_k$ ne s'applatissent pas lorsque $k\mapsto\infty$, alors
la dimension num\'erique de $K_X$ vaut z\'ero (la preuve de cette affirmation est tr\`es 
voisine de celle qui sera donn\'ee pour le th\'eor\`eme 2, donc on la laisse 
au soin du lecteur). 

\vskip 15pt
Le dernier paragraphe de cet article concerne
les applications quasi-conformes.
Par analogie avec la pseudo-forme volume de Kobayashi, on introduit
au troisi\`eme chapitre la notion de {\sl pseudo-forme volume quasi-conforme}. Comme cons\'equence du th\'eor\`eme
2 on obtient:
\medskip 
\noindent {\bf Corollaire. } {\sl Soit $X$ une vari\'et\'e k\"ahl\'erienne compacte, telle que 
$K_X$ soit nef. Si la pseudo-forme volume quasi-conforme est non-d\'eg\'en\'er\'ee, alors
$\nu (K_X)\geq 1$.}

\noindent Ce corollaire nous a \'et\'e
sugger\'e par Y.-T. Siu, que nous remercions vivement.

Dans le m\^eme esprit, une application $\varphi:\bC^n\mapsto X$ est dite
{\sl quasi-conforme en moyenne} s'il existe une constante positive $\delta$ telle que
$$\int_0^r{{dt}\over {t^{2n-1}}}\int_{B(t)}\vert J(\varphi)\vert _\omega^{2/n}d\lambda\geq 
\delta T_\omega(\varphi, r)$$

\noindent pour tout $r\gg 0$; on note $J(\varphi)$
le jacobien de l'application $\varphi$.
Pour une telle application, nous allons montrer le th\'eor\`eme suivant:

\medskip
\noindent {\bf Th\'eor\`eme 3.} {\sl Soit $(X, \omega)$ une vari\'et\'e k\"ahl\'erienne 
compacte de dimension $n$, telle que $K_X$ est pseudo-effectif. S'il existe une application 
holomorphe $\varphi: \bC^n\mapsto X$, quasi-conforme en moyenne, alors $\nu(K_X)= 0$.}

\noindent On remarquera que les th\'eor\`emes 2 et 3 etablissent des cas 
particuliers de la conjecture propos\'ee ci-dessus.
\vskip 3cm

\noindent {\twelvebf 1. Applications de degr\'e moyen \`a croissance lente}

\medskip
\noindent Consid\'erons une application holomorphe $\varphi: \Delta_r\rightarrow X$, o\`u
$X$ est une vari\'et\'e complexe compacte de dimension $n$, et $\Delta_r$ d\'esigne le polydisque 
de rayon $r$ dans $\bC^n$. On munit $X$ d'une m\'etrique hermitienne $\omega$; 
la quantit\'e: 

$$\deg(\varphi\vert \Delta_r):= \int_{\Delta_r}\varphi^*\omega^n/\int_X\omega^n$$

peut \^etre 
interpr\'et\'ee comme le degr\'e moyen de $\varphi$ sur ce polydisque.

\noindent Dans ce contexte, notre r\'esultat est le suivant:
\medskip
\noindent {\bf Th\'eor\`eme 1.} {\sl Soit $(X, \omega)$ une vari\'et\'e complexe
compacte de dimension $n$. S'il existe une application holomorphe 
$\varphi:\bC^n\rightarrow X$ non d\'eg\'en\'er\'ee,
et telle que:

$$\lim_{r\mapsto \infty}1/r^{2n}\deg(\varphi\vert \Delta_r)= 0,\leqno (1) $$

\noindent alors le fibr\'e canonique de $X$ n'est pas pseudo-effectif.}

\medskip

\noindent Remarquons que le produit d'un tore de dimension $n-1$ par $\bP^1$ montre que l'exposant $2n$ est  optimal comme entier pair. 

Comme on l'a d\'ej\`a rappel\'e dans l'introduction , Kodaira montre dans
[6] que sous les hypoth\`eses du th\'eor\`eme 1, on a $H^0(X, K_X^m)= 0$, pour tout
$m\geq 1$. Maintenant si $L$ est un fibr\'e en droites sur $X$, tel que 
pour un certain $m\geq 0$ on ait $H^0(X, L^m)\neq 0$, alors $L$ est pseudo-effectif.
Ainsi, le th\'eor\`eme 1 peut-etre vu comme g\'en\'eralisation du th\'eor\`eme
de Kodaira.
D'autre part, bien qu'en g\'en\'eral la pseudo-effectivit\'e 
d'un fibr\'e ne soit pas equivalente \`a
l'existence de sections non-nulles pour une de ses puissances, 
dans le cas du fibr\'e canonique
cela devrait \^etre vrai, d'apr\`es la conjecture d'abondance.

\

\noindent {\bf  Preuve.} Supposons que $K_X$ soit pseudo-effectif; alors il existe une fonction
$f\in L^1(X)$, semi-continue sup\'erieurement, telle que:

$$\Theta_h(K_X)+ i\ddbar f\geq 0\leqno (2)$$

\noindent au sens des courants sur la vari\'et\'e $X$, o\`u $h$ est la m\'etrique duale de
$\det (\omega)$ sur le fibr\'e canonique (voir [3], par exemple).

Prenons l'image inverse de (2) sur $\bC^n$ via l'application $\varphi$; on a

$$i\ddbar (\log   
(\vert \vert J_\omega(\varphi) \vert \vert ^2e^{f\circ \varphi}))\geq 0\leqno (3)$$

\noindent et alors la fonction $\tau:\bC^n\rightarrow \bR_+$ d\'efinie par
$\tau(z):=\vert \vert J(\varphi, z) \vert \vert _\omega^2e^{f\circ \varphi(z)} $
sera psh sur $\bC^n$. On rappelle maintenant la propri\'et\'e 
de convexit\'e suivante des fonctions psh (voir [11], par exemple).

\medskip
\noindent {\bf Lemme. }{\sl La fonction
$$M_\tau (t_1, ..., t_n):= \int _{[0, 2\pi]^n}\tau (t_1e^{i\theta_1},..., t_ne^{i\theta_n})d\theta$$

\noindent est croissante en chaque variable, et convexe en $\log (t_j)$.}

\medskip

\noindent Gr\^ace a ce lemme, on a la suite d'in\'egalit\'es 

$$\eqalign{
\int_{\Delta_r}\tau(z)d\lambda(z)= & \int_{[0, r]^n}t_1...t_nM_\tau (t_1, ..., t_n)dt_1...dt_n\geq \cr
\geq & \int_{[r_0, r]^n}t_1...t_nM_\tau (t_1, ..., t_n)dt_1...dt_n\geq \cr
\geq & M_\tau (r_0,..., r_0) \int_{[r_0, r]^n}t_1...t_ndt_1...dt_n=\cr
= & \bigl({{r^2- r_0^2}\over {2}}\bigr)^nM_\tau (r_0,..., r_0).\cr
}$$
Par ailleurs, la fonction $f$ est born\'ee sup\'erieurement sur $X$, et donc

$$\int_{\Delta_r}\tau(z)d\lambda(z)\leq C\int_{\Delta_r}
\vert \vert J_\omega(\varphi) \vert \vert _z^2d\lambda= \int_{\Delta_r}\varphi^*\omega^n$$

\noindent La condition de croissance impos\'ee par hypoth\`ese sur le degr\'e moyen
de l'applica\-ti\-on $\varphi$ montre que $M_\tau (r_0,..., r_0)= 0$, et ceci est vrai pour tout 
rayon $r_0> 0$. Notre hypoth\`ese permet \'egalement de changer l'origine de $\bC^n$, 
quitte \`a faire une translation, donc en resum\'e $\tau\equiv 0$. Mais ceci est 
clairement impossible, car l'image de l'application $\varphi$ contient un 
ouvert de la vari\'et\'e $X$. La contradiction ainsi obtenue montre que $K_X$ n'est pas pseudo-effectif,
et notre th\'eor\`eme est d\'emontr\'e.

\vskip 15pt

\noindent Supposons \`a pr\'esent que $X$ est une vari\'et\'e projective, de dimension $n$.
D'apr\`es les r\'esultats de [3], on sait que si le fibr\'e canonique n'est pas pseudo-effectif,
alors $X$ est uniregl\'ee. Le th\'eor\`eme 3 admet donc le corollaire suivant.
\medskip
\noindent {\bf Corollaire} {\sl Soit $X$ une vari\'et\'e projective
de dimension $n$. S'il existe une application holomorphe 
$\varphi:\bC^n\rightarrow X$ non d\'eg\'en\'er\'ee,
et telle que:

$$\lim_{r\mapsto \infty}1/r^{2n}\deg(\varphi\vert \Delta_r)= 0, $$

\noindent alors $X$ est unir\'egl\'ee.}

\medskip

\noindent Conjecturalement, toute compactification de $\bC^n$ devrait \^etre
rationnelle. Donc on peut voir le corollaire pr\'ec\'edent comme un premier
pas vers cette conjecture. Malheureusement, la fa\c con dont les courbes rationnelles
apparaissent dans notre r\'esultat (i.e. via la th\'eorie de Mori) n'est pas du tout explicite.

Il est d'ailleurs possible qu'aucune courbe alg\'ebrique de $\bC^n$ n'ait pour image une courbe alg\'ebrique de $X$, m\^eme si $\deg(\varphi\vert \Delta_r)$ est uniform\'ement born\'ee: soit $\Omega\cong \bC^n\subsetneq \bC^n\subset \bP^n=X$ un domaine de Fatou-Bieberbach obtenu comme bassin d'attraction d'un automorphisme r\'egulier (voir [12], 2.2, p. 124, et p.125 pour des exemples) de $\bC^n,n\geq 3$. Alors $\Omega$ ne contient aucune courbe alg\'ebrique (par [12], Th\'eor\`eme 2.4.4 et remarque 2.4.5 (1) , p. 136), et on est dans la situation \'evoqu\'ee ci-dessus. 

La m\'ethode de d\'emonstration du th\'eor\`eme 1 pr\'ec\'edent devrait permettre d' \'etablir la g\'en\'eralisation suivante (illustr\'ee par le produit $X=\bP^{n-p+1}\times A$, o\`u $A$ est une vari\'et\'e ab\'elienne de dimension $p-1$):

{\bf Question:} {\sl Soit $X$ une vari\'et\'e projective
de dimension $n$. S'il existe une application holomorphe 
$\varphi:\bC^n\rightarrow X$ non d\'eg\'en\'er\'ee,
et telle que:

$$\lim_{r\mapsto \infty}1/r^{2p}\deg(\varphi\vert \Delta_r)= 0, $$

\noindent et si $f:X\mapsto Y$ est une application m\'eromorphe surjective telle que $dim(Y)=p$, alors:  $Y$ est-elle unir\'egl\'ee?}

Une r\'eponse affirmative \`a cette question montrerait que le {\it quotient rationnel} de $X$ (voir [4] et [10]) est de dimension au plus $(p-1)$. En particulier, si $\deg(\varphi\vert \Delta_r)$ est born\'e, $X$ (et donc toute compactification de $\bC^n$) serait du moins rationnellement connexe.

\vskip 30pt

\noindent {\twelvebf 2. Preuve du th\'eor\`eme 2}

\vskip 15pt  {\bf 2.1} On commence par rappeler la notion de dimension num\'erique d'une 
classe de cohomologie pseudo-effective, telle qu'elle a \'et\'e introduite par 
Boucksom, Tsuji dans [1], [2], [13].

Soit $\{\alpha\}\in H^{1, 1}(X, \bR)$ une classe de cohomologie pseudo-effective
(i.e., il existe un courant positif ferm\'e $T\in \{\alpha\}$, voir [3] pour une 
pr\'esentation plus compl\`ete de cette notion). Pour chaque nombre r\'eel 
$\varepsilon> 0$, la classe de cohomologie $\{\alpha+ \varepsilon \omega\}$
est ``big'', et il existe un courant positif ferm\'e $T_\varepsilon\in 
\{\alpha+ \varepsilon \omega\}$, tel que ses singularit\'es sont concentr\'ees le long
d'un ensemble analytique $Y_\varepsilon$. On d\'efinit le {\sl nombre 
d'intersection mobile} (voir [2]) comme suit

$$(\alpha^k\wedge \omega^{n-k}):= \lim_{\varepsilon\mapsto 0}\sup_ 
{T_\varepsilon\in 
\{\alpha+ \varepsilon \omega\}}\int_{X\setminus Y_\varepsilon}
T_\varepsilon^k\wedge \omega^{n-k}.$$

\medskip
\noindent {\bf D\'efinition {\rm ([2])}} {\sl Soit 
 $\{\alpha\}\in H^{1, 1}(X, \bR)$ une classe de cohomologie pseudo-effective sur
$X$. La dimension num\'erique de $\alpha$ est d\'efinie par 
$$\nu(\alpha):= \max \{k\in \bZ/(\alpha^k\wedge \omega^{n- k}\neq 0\}.$$}
\medskip

Il a \'et\'e demontr\'e par S.Boucksom (voir [1]) qu'une classe
$\{\alpha\}$ est de dimension num\'erique maximale si et seulement si
elle contient un r\'epr\'esentant strictement positif. 
Concernant l'autre cas extr\^eme $\nu(\alpha)= 0$, la situation est malheureusement
beaucoup moins bien comprise (cf. [1] pour quelques r\'esultats dans cette direction).

\

{\bf 2.2} Supposons \`a pr\'esent qu'on ait
$\nu:= \nu(K_X)\geq 1$. Alors pour tout $\varepsilon> 0$, on a un courant positif
$T_\varepsilon\in c_1(K_X)+ \varepsilon \{\omega\}$, dont les singularit\'es sont
concentr\'ees le long d'un ensemble analytique $Y_\varepsilon$, tel que 
$$\int_{X\setminus Y_\ep}T_\varepsilon^\nu \wedge \omega^{n- \nu}\geq \delta_0> 0$$

\noindent uniform\'ement par rapport \`a $\varepsilon$. Autrement dit, il existe une
famille de modifications $\mu_\varepsilon: X_\varepsilon\to X$ telle que 
$\mu_\varepsilon^*T_\varepsilon= [E_\varepsilon]+ \widetilde\alpha_\varepsilon$, o\`u
$E_\varepsilon$ est un $\bQ$-diviseur effectif, (dont la partie $\mu_{\varepsilon}$-exceptionnelle
provient des singularit\'es de $T_\ep$ en codimension 2), 
et $\widetilde\alpha_\varepsilon$
est une $(1, 1)$--forme semi-positive sur $X_\varepsilon$, tels que
$$\int_{X_\varepsilon}\widetilde\alpha_\varepsilon^\nu \wedge\mu_\varepsilon^*\omega^{n- \nu}\geq 
\delta_0> 0.$$

\

\noindent {\bf 2.3} On rappelle maintenant que si $\mu: X_1\rightarrow X$ est un 
\'eclatement de centre lisse $Y\subset X$, alors il existe une m\'etrique $h$ sur
le fibr\'e associ\'e au diviseur exceptionnel $E$ telle que pour tout 
$0<\delta \ll 1$, la forme diff\'erentielle $\wt\omega_\ep:=\mu^*\omega- \delta \Theta_h(E)$
est positive d\'efinie sur $X_1$.   
Ainsi, pour chaque $\varepsilon> 0$, on construit $\wt\omega_\ep$ une 
m\'etrique k\"ahl\'erienne sur $X_\ep$, telle que:

\item {(i)} $\wt\omega_\ep- \mu_\ep^*\omega$
soit un multiple (n\'egatif) du diviseur exceptionnel.
\smallskip
\item {(ii)} en chaque point de $X_\ep$, on a $\wt\omega_\ep\geq 1/2\mu_\ep^*\omega$. 
\smallskip
\item {(iii)} le volume de $(X_\ep, \wt\omega_\ep)$
soit major\'e par le volume de $(X, 2^{1/n}\omega)$.

\noindent En effet, les conditions (i)-(iii) sont clairement satisfaites dans le
cas d'un seul \'eclatement de centre lisse (quitte \`a choisir le param\`etre 
$\delta$ ci-dessus assez petit). Le cas g\'en\'eral est d\'eduit du fait qu'on peut
supposer que la modification $\mu_\ep$ est une compos\'ee d'\'eclatements de centres lisses.

\  

\noindent {\bf 2.4} Nous utilisons maintenant le th\'eor\`eme de Yau [14], afin d'obtenir dans la classe 
$\{\wt\alpha_\ep+ \ep\wt\omega_\ep\}$ un r\'epr\'esentant dont le d\'eterminant est 
constant par rapport \`a $\wt\omega_\ep$. 
\medskip
\noindent {\bf Th\'eor\`eme {\rm (Yau).}} {\sl Soit $(X, \beta)$ une vari\'et\'e k\"ahl\'erienne
compacte, de dimension $n$, et soit $dV$ un \'el\'ement volume sur $X$, tel que
$\int_X\beta^n= \int_XdV$. Alors il existe $\rho\in {\cal C}^\infty(X)$
telle que
\item {(i)} $\beta+ i\ddbar \rho > 0$ sur $X$.
\smallskip
\item {(ii)} $(\beta+ i\ddbar \rho)^n= dV$.}
\medskip
\noindent On applique le th\'eor\`eme
pr\'ec\'edent pour chaque $X_\ep$, munie de la m\'etrique k\"a\-hl\'e\-rienne
$\wt\alpha_\ep+ \ep\wt\omega_\ep$, et l'\'el\'ement volume $dV:= C(\ep)\wt\omega_\ep^{n}$.
pour une certaine constante ad\'equate de normalisation $C(\ep)$, dont le role est 
de satisfaire la condition cohomologique du th\'eor\`eme de Yau.
Ainsi on montre l'existence d'une famille de fonctions
$\rho_\ep\in {\cal C}^\infty(X_\ep)$, telle que

\item {(a)} $\wt\alpha_\ep+ \ep\wt\omega_\ep + i\ddbar\rho_\ep> 0$ sur $X_\ep$.
\smallskip
\item {(b)}  $(\wt\alpha_\ep+ \ep\wt\omega_\ep + i\ddbar\rho_\ep)^n= 
C(\ep)\wt\omega_\ep^{n}$.

\noindent Nous observons maintenant que, gr\^ace au th\'eor\`eme de Yau, les
propri\'et\'es nu\-m\'e\-riques de la classe canonique se refl\`etent dans ses
propri\'et\'es m\'etriques, i.e. on obtient un minorant pour la constante $C(\ep)$
en int\'egrant l'\'egalit\'e (b) ci-dessus. Ainsi, on montre l'existence d'une constante 
$C_1> 0$, telle que $C(\varepsilon)\geq C_1\ep^{n- \nu}$, lorsque $\ep\to 0$. En effet,
l'\'egalit\'e (b) implique

$$
\eqalign{
C(\ep)\Vol(X_\ep, \wt\omega_\ep)= & \int_{X_\ep}\bigl(\wt\alpha_\ep+ \ep\wt\omega_\ep\bigr)^n\geq\cr
\geq & C_n^\nu \ep^{n-\nu}\int_{X_\ep}\wt\alpha_\ep^\nu\wedge \wt\omega_\ep^{n-\nu}\geq C_n^\nu/2^{n-\nu}
\ep^{n-\nu}\int_{X_\ep}\wt\alpha_\ep^\nu\wedge \mu_\ep^*\omega^{n-\nu}\geq\cr
\geq & C_0\ep^{n-\nu}\cr
}$$
(dans la suite des in\'egalit\'es pr\'ec\'edentes, on a utilis\'e la semi-positivit\'e de 
$\alpha_\ep$, ainsi que la d\'efinition de la dimension num\'erique dans le cadre
pseudo-effectif). L'existence de $C_1$ se d\'eduit du fait que le volume de 
$(X_\ep, \wt\omega_\ep)$ est uniform\'ement major\'e.

\

\noindent {\bf 2.5} Consid\'erons le courant image directe
 $\Theta_\ep=\mu_{\ep,*}(\wt\alpha_\ep+ \ep \wt\omega_\ep + [E_\ep]+ i \ddbar\rho_\ep)$
sur la vari\'et\'e $X$. C'est un courant positif ferm\'e et sa classe de cohomologie est
$c_1(K_X)+ 2\ep\{\omega\}$, comme on le voit imm\'ediatement. De plus, $\Theta_\ep$ est 
non-singulier sur $X\backslash Y_\ep$ et l'\'equation de Calabi-Yau (ii) qu'on r\'esout sur
$X_\ep$ montre l'existence d'une constante $C_2> 0$ telle que

$$\Theta_\ep^n\geq C_2\cdot\ep^{n-\nu}\cdot\omega^n\leqno (4)$$

\noindent en tout point de $X\backslash Y_\ep$. Pour v\'erifier la relation 
(4), pla\c cons-nous en un point $x_0\in X\setminus Y_\ep$. Au voisinage de 
$x_0$, on obtient 
$\Theta_\ep=\mu_{\ep,*}(\wt\alpha_\ep+ \ep \wt\omega_\ep + i \ddbar\rho_\ep)$, car 
le support de $[E_\ep]$ est contenu dans l'image inverse de $Y_\ep$. L'\'egalit\'e 
(b) combin\'ee avec la propri\'et\'e (ii) de la famille de m\'etriques $(\omega_\ep)$ 
et le minorant de $C(\ep)$ montrent que 

$$(\wt\alpha_\ep+ \ep \wt\omega_\ep + i \ddbar\rho_\ep)^n
\geq C_1\ep^{n-\nu}/2^n\mu_\ep^*\omega^n.\leqno (5)$$
Maintenant au voisinage de $x_0$, l'application $\mu_\ep$ est un isomorphisme, et
ainsi l'in\'egalit\'e (5) implique (4), avec une constante $C_2:= C_1/2^n$.

\

{\bf 2.6} On \'ecrit maintenant $\Theta_\ep=\alpha + 2 \ep\omega + i\ddbar f_\ep$, pour une certaine
fonction $f_\ep \in L^1(X)\cap \cC^\infty (X\backslash Y_\ep)$, normalis\'ee de telle sorte que
$\displaystyle\int_X f_\ep \omega^n = 0$ (pour simplifier l'\'ecriture, on a not\'e
$\alpha:= \Theta_{\det(\omega)}(K_X)$). Le but de la normalisation de $f_\ep$ 
est d'obtenir la relation d'uniformit\'e suivante.

\medskip

\noindent {\bf Lemme.} {\sl Il existe une constante positive $C_3:= C_3(X, \alpha, \omega) $ 
telle que

$$\max \Bigl(\int_X\vert f_\ep\vert dV_\omega,  \sup_X(f_\ep)\Bigr)\leq C_3\leqno (6)$$

\noindent pour tout $\ep> 0$.}

\medskip 
\noindent {\sl Preuve du lemme.} Tout d'abord, le lemme est bien connu dans la cadre suivant:
soit $\beta$ une (1, 1)-forme diff\'erentielle de classe $\cC^\infty$ sur $X$, et $f\in \cC^\infty(X)$
une fonction telle que $\beta+ i\ddbar f\geq 0$, et telle que $\int_XfdV_\omega= 0$. Alors 
 
$$\max \Bigl(\int_X\vert f\vert dV_\omega,  \sup_X(f)\Bigr)\leq C(X, \beta, \omega)\leqno (7)$$
(voir par exemple [5]). L'\'enonc\'e g\'en\'eral est con\-s\'e\-quence
de ce fait, et du th\'eor\`eme de regularisation suivant, d\^u \`a Demailly:
\medskip 
\noindent {\bf Th\'eor\`eme {(\rm ([6]).)}} {\sl Soit $X$ une vari\'et\'e 
complexe compacte et soit $T= \alpha+ i\ddbar f$ un courant positif ferm\'e de type $(1,1)$.
Alors pour chaque entier positif $k$, il existe une fonction 
$f_k\in \cC^\infty(X)$ telle que:

\item {(1)} $f_k\to f$ en norme $L^1$ et ponctuellement sur $X$. Posant, de plus: $T_k= \alpha+ i\ddbar f_k$, alors: 
\smallskip

\item {(2)} $T_k\in \{\alpha \}$, et $T_k\geq -\lambda_k\omega$, o\`u $\lambda_k(x)\to \nu(T, x)$ lorsque $k\to\infty$ (autrement dit, en chaque point $x\in X$, la perte de positivit\'e est de la taille du
nombre de Lelong $\nu(T, x)$ du courant initial $T$).
}
\medskip 

\noindent En effet, on applique \`a chaque $f_\ep$ le th\'eor\`eme pr\'ec\'edent;
l'observation importante est la suivante. Etant donn\'e que la classe de cohomologie
des courants $\Theta_\ep$ est born\'ee par rapport \`a $\ep$, les nombres de Lelong de 
$\Theta_\ep$ le sont \'egalement, car ils sont domin\'es par la masse du courant,
et dans le cadre k\"ahl\'erien, cette quantit\'e est cohomologique 
(ceci reste, en fait, vrai pour les courants de type $(1,1)$ sur une vari\'et\'e complexe compacte, 
gr\^ace \`a l'existence des m\'etriques de Gauduchon). Donc, la perte de positivit\'e pour les courants
r\'egularis\'es est uniforme par rapport aux param\`etres $\ep, k$, et ainsi notre lemme 
est d\'emontr\'e par le r\'esultat rappel\'e au d\'ebut de la preuve.
\vskip 15pt

\noindent {\bf 2.7} Revenons maintenant \`a notre application $\varphi : \bC^n \to X$.
Si on note $\Delta_1$ le polydisque unit\'e, alors 
$B_\omega(x_0, \delta)\subset \varphi(\Delta_1)$ pour un certain rayon $\delta> 0$,
et la relation (1) montre l'existence d'un $x_\ep\in \Delta_1$ tel que, pour tout $\ep> 0$:

$$  f_\ep\circ \varphi(x_\ep)\geq -C'_3/\delta^{2n}\leqno (7')$$

Afin de ne pas trop alourdir les notations, on va supposer que
$x_\ep= 0$, l'origine de $\bC^n$ (les arguments present\'es par la suite
montreront qu'on peut se le permettre, car translater l'origine de $\bC^n$ en $x_\ep
$ revient \`a changer $r$ en $ r+1$ \`a la fin de la preuve, et la conclusion d\'esir\'ee
n'est pas affect\'ee par ce changement). 

Sur $\bC^n$, consid\'erons le jacobien $J(\varphi)\in H^0(\bC^n, \varphi^*K_X^{-1})$. 
La norme de $J(\varphi)$ par rapport \`a la m\'etrique $\det (\omega)$ sur $K_X^{-1}$ est donn\'ee 
par l'\'egalit\'e $\displaystyle \varphi^*\omega^n= 
\vert \vert J(\varphi)\vert \vert _\omega^2d\lambda$,
o\`u $d\lambda$ est la mesure de Lebesgue de $\bC^n$. Ainsi, on a 
$$i\ddbar \log \Bigl(\Vert J(\varphi)\Vert ^2_{\omega}
e^{f_\ep\circ\varphi}\Bigr)+ 2\ep\varphi^*\omega\geq \varphi^*\Theta_\ep$$

\noindent car la diff\'erence est donn\'ee par la courant d'int\'egration sur 
le lieu des z\'eros du jacobien. 

Maintenant, en tout  
 point $z\in \bC^n\backslash \varphi^{-1}(Y_\ep)$, on a l'in\'egalit\'e
$$\varphi^*\Theta_\ep \wedge (i\ddbar\Vert z\Vert^2)^{n-1}\geq C_2^{1/n}\ep^{1- {{\nu}\over{n}}}
 \Vert J(\varphi)\Vert ^{2\over n}\cdot (i\ddbar\Vert z\Vert^2)^{n} $$
Pour v\'erifier cette affirmation, soient $\lambda_\ep^{(j)}(z)$ les valeurs propres de 
$\varphi^*\Theta_\ep$ par rapport \`a la m\'etrique euclidienne au point $z\in \bC$. 
La relation (4) montre que 
$$\prod\lambda_\ep^{(j)}(z)\geq C_2\ep^{n-\nu}\vert \vert J(\varphi, z)\vert \vert _\omega^{2}$$
et par l'in\'egalit\'e de la moyenne on d\'eduit un minorant pour la somme des valeurs propres, qui est 
pr\'ecisement l'in\'egalit\'e ci-dessus.
De plus, on remarque que cette relation reste vraie sur $\bC^n$ tout entier, au sens
des distributions. En conclusion on obtient
$$\big (i\ddbar \log \Vert J(\varphi)\Vert ^2_{\omega}
        e^{f_\ep\circ\varphi}+2\ep\varphi^*\omega\big ) \wedge 
(i\ddbar\Vert z\Vert^2)^{n-1} \geq \ep^{1- {{\nu}\over{n}}}
 \Vert J(\varphi)\Vert ^{2\over n}_{\omega}d\lambda\leqno (8)$$

\

\noindent {\bf 2.8} Pour la suite, on voudrait utiliser l'in\'egalit\'e (8)
afin d'appliquer les arguments de courbure n\'egative (cf. Ahlfors, Kodaira, Griffiths)
mais malheureusement, on ne peut pas le faire directement, \`a cause du terme
$\varphi^*\omega \wedge (i\ddbar \Vert z\Vert ^2)^{n-1}$ (la trace de la m\'etrique image inverse 
par rapport \`a la m\'etrique euclidienne). Afin d'inclure ce terme dans la deriv\'ee 
logarithmique ci-dessus, on r\'esout le probl\`eme de Dirichlet suivant.

$$
\left\lbrace
\eqalign{
& \Delta\Psi_r= (\varphi^*\omega)\wedge (i\ddbar\Vert z\Vert^2)^{n-1}/d\lambda
  \ \ dans\ \ \Vert z\Vert \leq r \cr
& \Psi_r=0  \ \ sur\ \ \Vert z\Vert = r.\cr
}
\right.
$$
La solution $\Psi_r$ est d\'etermin\'ee de mani\`ere unique par la formule de Green; 
on laisse au soin du lecteur de v\'erifier que la valeur au bord qu'on a consid\'er\'ee est 
optimale pour ce qui va suivre.

A l'aide de la fonction $\Psi_r$ l'in\'egalit\'e (8) s'\'ecrit 
$$
\Delta\left ( \log \Vert J(\varphi)\Vert ^{2\over n}_{\omega}
e^{{\ep\Psi_r + f_\ep\circ\varphi}\over {n}}
\right )\geq {{C_2^{1/n}}\over {n}}\ep^{1-{\nu\over n}}\Vert J(\varphi)\Vert _\omega ^{2\over n}
\leqno (9)$$
Gr\^ace \`a l'in\'egalit\'e (6), la derni\`ere relation implique

$$
\Delta\left ( \log \Vert J(\varphi)\Vert ^{2\over n}_{\omega}
e^{{\ep\Psi_r + f_\ep\circ\varphi}\over {n}}\right )
\geq C_4\ep^{1-{\nu\over n}}\Vert J(\varphi)\Vert _\omega ^{2\over n}e^{{f_\ep\circ\pi}\over {n}}
\leqno (10)$$
(l'expression de la constante $C_4$ se d\'eduit facilement de $C_2$ et $C_3$).

Par ailleurs, la fonction $\Psi_r$ solution du probl\`eme de Dirichlet 
est sous-harmo\-nique, donc par le principe du maximum,
$\displaystyle\max_{\Vert z\Vert \leq r}\Psi_r(z)=0$, et finalement on peut \'ecrire
$$\Delta\left ( \log \Vert J(\varphi)\Vert ^{2\over n}_{\omega}
e^{{\ep\Psi_r + f_\ep\circ\varphi}\over {n}}
\right )\geq  C\ep^{1-{\nu\over n}}\Vert J(\varphi)\Vert ^{2\over n}_{\omega}
e^{{\ep\Psi_r + f_\ep\circ\varphi}\over {n}}
\leqno (11)$$
Les arguments de courbure n\'egative  suivants sont classiques (voir par exemple [5], [6]).
\noindent
Sur $\Vert z\Vert \leq r$, consid\'erons la fonction $\tau_\ep$ d\'efinie par :
$$z\to \left (1-{{\Vert z\Vert^2}\over{r^2}} \right )^2
\Vert J(\varphi, z)\Vert _\omega ^{2\over n}
e^{{\ep\Psi_r(z) + f_\ep\circ\varphi(z)}\over {n}}:=\tau_\ep (z)$$
C'est une fonction positive, et sur la sph\`ere $\left ( \Vert z \Vert=r \right )$ 
elle vaut z\'ero; par cons\'e\-quent, son point de maximum est atteint en un $z_{max}$ de $\Vert z\Vert < r$.
\noindent Alors on a
$$
\Delta \log \tau_\ep (z_{max})\leq 0,
$$ 
Par ailleurs, un calcul sans difficult\'e montre que: 
$$i\ddbar\log \left (1-{{\Vert z\Vert^2}\over{r^2}} \right ) 
\wedge \left (i\ddbar\Vert z\Vert^2 \right )^{n-1}\geq
-{{d\lambda}\over{r^2(1-{{\Vert z\Vert^2}\over{r^2}})^2}}
\leqno (12)$$
Donc au point de maximum $z= z_{max}$ on aura:
$$
C_4\ep^{1-{\nu\over n}}\Vert J(\varphi, z)\Vert ^{2\over n}_{\omega}
e^{{\ep\Psi_r(z) + f_\ep\circ\varphi(z)}\over {n}}\leq
{{1}\over{r^2(1-{{\Vert z\Vert^2}\over{r^2}})^2}},
$$
comme cons\'equence de (11) et (12). Alors on obtient:
$$
C_4\ep^{1-{\nu\over n}}\Vert J(\varphi, z)\Vert ^{2\over n}_{\omega}
e^{{\ep\Psi_r(z) + f_\ep\circ\varphi(z)}\over {n}}\cdot
(1-{{\Vert z\Vert^2}\over{r^2}})^2\leq
{{1}\over{r^2}}
$$
au point $z_{max}$, et donc $\displaystyle \tau_\ep(z)\leq {{1}\over{C_4\ep^{1-{\nu\over n}}r^2}}$ pour $z\in\bB(r)$.
En particulier, si $z=0$,  l'in\'egalit\'e pr\'ec\'edente devient:
$$\displaystyle\Vert J(\varphi, 0)\Vert ^{2\over n}_{\omega}\ep^{1-{\nu\over n}}
e^{{\ep\Psi_r(0) + f_\ep\circ \varphi(0)}\over {n}}\leq
{{1}\over{r^2}}.$$

\noindent Afin de d\'eterminer la quantit\'e $\Psi_r(0)$, on applique la formule de Green

$$\int_{\Vert z \Vert= r}\Psi_r(\xi)d\sigma- \Psi_r(0)= 
\int_0^r{{dt}\over {t^{2n-1}}}\int _{B(t)}\Delta \Psi_r d\lambda$$

\noindent et comme $\Psi_r$ est la solution du probl\`eme de Dirichlet,
on voit que $\Psi_r(0)=-T(\varphi,r)$, l'indicatrice de croissance de 
la fonction $\varphi$. En conclusion, on a 
$$
\Vert J(\varphi, 0)\Vert _\omega ^{2\over n}\leq{{C_5}\over{r^2}}\cdot{{1}\over{\ep^{1-{\nu\over n}}}}
\cdot e^{\ep T(\varphi,r)},
$$
compte-tenu de l'in\'egalit\'e (7').

\noindent On peut remarquer que dans l'expression pr\'ec\'edente,
la constante $C_5$ d\'epend uniquement des quantit\'es suivantes: g\'eom\'etrie de $(X, \omega)$, et
rayon de la boule euclidienne contenue dans l'image $\varphi (\Delta_1)$.
 
Maintenant on observe que les param\`etres $\ep$ et 
$r$ sont  ind\'e\-pendants. Par cons\'equent, on peut prendre
$\ep:= T^{-1}(\varphi, r)$ et on obtient

$$\Vert J(\varphi, 0)\Vert_\omega ^{2\over n}\leq{{C}\over{r^2}}\cdot(T(\varphi,r))^{1-{\nu\over n}}$$
La preuve du th\'eor\`eme 1 est ainsi compl\`ete.

\medskip
\noindent{\bf Remarque.} Comme cons\'equence de la m\'ethode de d\'emonstration pr\'ec\'edente, 
on obtient l'enonc\'e suivant:
\smallskip
\noindent{\bf Th\'eor\`eme $2^\prime $}
Soit $(X^n,\omega)$ une vari\'et\'e k\"ahl\'erienne compacte de dimension $n$, dont le fibr\'e canonique
$K_X$ est est pseudo-effectif, avec $\nu:=\nu(K_X)$. 

Si $\varphi : B^n(r)\to X $ est d\'efinie sur une boule de rayon $r>1$ de $\Bbb C$$^n$, et si $\varphi(B(1))\supset B_\omega(x_0, \delta_0)$ pour un certain $x_0\in X$, alors 
$$
T(\varphi, r^\prime)^{1-\nu/n}\geq C(X,\omega,\delta_0)\cdot (r')^2\cdot 
\Vert J(\varphi, 0)\Vert ^{2\over n}_{\omega} \ \ \forall\ \ r^\prime\leq r.
$$

\medskip
\noindent Pour faire le lien avec la premi\`ere partie de cet article, 
il serait int\'eressant d'ana\-lyser les cons\'equences de la positivit\'e num\'erique de 
$K_X$ sur la croissance du degr\'e moyen de l'application $\varphi$
(i.e., d'obtenir l'analogue du th\'eor\`eme 2 pour le degr\'e moyen \`a la place de
la fonction caract\'eristique de $\varphi$).

\vskip 30pt
\noindent {\twelvebf 3. Applications quasi-conformes.}

\medskip
\noindent Nous rappelons qu'\'etant donn\'e  un point $x_0\in X$,
on peut consid\'erer la quantit\'e :
$$r(x_0)=sup\left\lbrace r>0\ \ /\ \ \exists\ \ f:\bB^n(r)\to X \ tel \ que\
\left\lbrace 
\eqalign{
& f(0)=x_0 \cr
&\Vert J(f, 0)\Vert_\omega=1\cr
}\right.
\right\rbrace$$

Alors le th\'eor\`eme de Kobayashi-Ochiai peut \^etre reformul\'e comme suit.
\smallskip
\noindent{\bf Th\'eor\`eme {\rm (KO)}} {\sl Soit $X$ une vari\'et\'e
projective, dont le fibr\'e canonique est big. Alors
$r(x_0)<\infty,\ \forall x_0\in X$.}
\smallskip
Comme cela nous a \'et\'e sug\'er\'e par Y-T Siu, on peut consid\'erer
une quantit\'e analogue o\`u seules les applications quasi-conformes
seront impliqu\'ees. On d\'efinit
$$r_{q,C}(x_0)=sup\left\lbrace r>0/ \exists\ \ f:\bB^n(r)\to X \ tel \ que\
\left\lbrace 
\eqalign{
& f(0)=x_0 \cr
&\Vert J(f, 0)\Vert _\omega=1\cr
&\Vert df(z)\Vert _\omega<C\cdot\Vert J(f, z)\Vert _\omega ^{1\over n}
}\right.
\right\rbrace$$
Dans ce contexte, la m\'ethode de d\'emonstration du th\'eor\`eme 2 donne le r\'esultat suivant.
\smallskip
\noindent{\bf Corollaire 3.1} {\sl Soit $(X,\omega)$ une vari\'et\'e k\"ahl\'erienne
compacte avec $K_X$ nef. Si $r_{q,C}(x_0)=+\infty$ pour un certain
$x_0\in X $, alors $\nu(K_X)=0$.}

\smallskip
\noindent{\sl Preuve. } L'hypoth\`ese implique 
l'existence d'une famille de boules $\varphi_r: \Delta_r\to X$ telle que
$$
\left\lbrace 
\eqalign{
&(1) \varphi_r (0)=x_0 \cr
&(2) \Vert J(\varphi_r, 0)\Vert= 1\cr
&(3) \Vert d\varphi_r(z)\Vert <C\cdot\Vert J(\varphi_r, z)\Vert^{1\over n}
}\right. $$

\noindent o\`u la derni\`ere condition est v\'erifi\'ee en tout point $z$ du polydisque
$\Delta_r$.

Le fibr\'e canonique de $X$ \'etant nef,  le th\'eor\`eme de 
Yau (voir la preuve du th\'eor\`eme 2) montre l'existence d'une suite de 
m\'etriques $h_\varepsilon:= he^{-f_\varepsilon}$, telles que
$f_\varepsilon\in{\cal C}^\infty (X)$, et telles que
$$(\Theta_{h}(K_X)+ i\ddbar f_\varepsilon+ \varepsilon\omega)^n= C_\varepsilon \omega^n
\leqno (13)$$

\noindent o\`u la constante de normalisation $C_\varepsilon$ est de l'ordre de  
${\cal O}(\varepsilon^{n-\nu})$, et o\`u $\nu= \nu(K_X)$ d\'esigne la dimension num\'erique du
fibr\'e canonique de $X$. 

Pour chaque $r> 0$, l'in\'egalit\'e de la moyenne combin\'ee avec (13)
implique

$$(i\ddbar \log \vert \vert J(\varphi_r) \vert \vert ^2e^{f_\varepsilon\circ\varphi_r}+
\varepsilon \varphi_r^*\omega )\wedge (i\ddbar \vert \vert z \vert \vert ^2)^{n-1}\geq
\varepsilon ^{1-\nu/n}\vert \vert J(\varphi_r) \vert \vert ^{2/n}d\lambda\leqno (14)$$

\noindent Supposons \`a pr\'esent que $\nu(K_X)\geq 1$; on aura alors clairement $\varepsilon^{1-\nu/n}\gg \varepsilon$,
lorsque $\varepsilon\mapsto 0$.
L'in\'egalit\'e
(14) et l'hypoth\`ese sur $\varphi_r$ donnent:

$$\Delta\log \vert \vert J(\varphi_r) \vert \vert ^{2/n}e^{f_\varepsilon \circ \varphi_r/n}
\geq \varepsilon^{1-\nu/n}\vert \vert J(\varphi_r) \vert \vert ^{2/n}e^{f_\varepsilon \circ\varphi_r/n}.
\leqno (15)$$
(car le terme $(\varphi_r^*\omega )\wedge (i\ddbar \vert \vert z \vert \vert ^2)^{n-1}$
est domin\'e --\`a une constante universelle pr\`es--par $\Vert J(\varphi_r)\Vert^{2\over n}d\lambda$, gr\^ace 
\`a la relation de quasi-conformalit\'e (3)).

\noindent
Fixons par la suite une valeur de $\varepsilon\ll 1$ telle que $(15)$ soit satisfaite. Le reste
de la preuve suit les arguments de courbure n\'egative dej\`a expos\'es en 2.9, et on obtient ainsi 
$$\vert \vert J(\varphi_r, 0)\vert \vert _\omega^{2/n}
e^{f_\varepsilon \circ\varphi_r(0)/n}\leq C/r^2.\leqno (16)$$

\noindent Maintenant les conditions (1) et (2) montrent que
la relation (16) est \'equivalente \`a l'in\'egalit\'e: $\displaystyle e^{f_\varepsilon(x_0)}\leq C/r^2$, et on obtient une contradiction  lorsque $r\mapsto \infty$. La preuve du corollaire est achev\'ee.

\vskip 15pt

\

Dans le cours de ce dernier paragraphe, nous allons \'etendre le r\'esultat 
pr\'e\-c\'edent au cas o\`u $\varphi$ est seulement suppos\'ee quasi-conforme {\it en moyenne}.

\

Soit $\varphi:\bC^n\rightarrow X$ une application holomorphe non-d\'eg\'en\'er\'ee.
On suppose l'existence d'une constante $\delta> 0$, telle que 
$$\int_0^r{{dt}\over {t^{2n-1}}}\int_{B(t)}\vert \vert J(\varphi) \vert \vert_{\omega} ^{2/n}
d\lambda\geq CT_\omega (\varphi, r)\leqno (17)$$

\noindent pour tout $r\gg 0$. Par exemple, si l'application $\varphi$ est quasi-conforme,
la condition $(17)$ est automatiquement satisfaite; en g\'en\'eral, appelons une telle 
application {\sl quasi-conforme en moyenne}. Nous allons maintenant pr\'esenter
une preuve de l'enonc\'e suivant.

\medskip
\noindent {\bf Th\'eor\`eme 3.2.} {\sl Soit $X$ une vari\'et\'e k\"ahl\'erienne
compacte de dimension $n$, et $\varphi: \bC^n\rightarrow X$ une application holomorphe
quasi-conforme en moyenne. Si le fibr\'e canonique de 
$X$ est nef, alors $\nu(K_X)= 0$.}

\medskip

\noindent{\bf Preuve.} Dans ce qui va suivre, on utilise de fa\c con essentielle
des arguments de la th\'eorie de Nevanlinna, et \'egalement les m\'etriques sur
le fibr\'e canonique qu'on a construites au cours de la preuve du th\'eor\`eme 2. 
Le point de d\'epart est l'in\'egalit\'e:

$$i\ddbar \log \vert \vert J(\varphi) \vert \vert ^2e^{f_\varepsilon \circ \varphi}+
\varepsilon \varphi^*\omega )\wedge (i\ddbar \vert \vert z \vert \vert ^2)^{n-1}\geq
\varepsilon ^{1-\nu/n}\vert \vert J(\varphi) \vert \vert ^{2/n}d\lambda\leqno (18)$$

\noindent en tout point de $\bC^n$. On int\`egre cette relation \`a la mani\`ere de
Nevanlinna; ainsi, pour tout $r> 0$ on aura:

$$\eqalign{
\int_0^r{{dt}\over {t^{2n-1}}}\int _{B(t)}
\Delta \log \vert \vert J(\varphi) \vert \vert_\omega ^2e^{f_\varepsilon \circ \varphi}d\lambda}$$

$$\geq  
\varepsilon ^{1-\nu/n}
\int_0^r{{dt}\over {t^{2n-1}}}\int _{B(t)}\vert \vert J(\varphi) \vert \vert _\omega ^{2/n}d\lambda 
 -\varepsilon T_\omega (\varphi, r).
$$

\noindent Si $\nu(K_X)> 0$, 
alors quitte \`a prendre $\varepsilon\ll 1$, (qui sera fix\'e par la suite), notre 
hypoth\`ese de quasi-conformalit\'e en moyenne et la relation ci-dessus entrainent:

$$\int_0^r{{dt}\over {t^{2n-1}}}\int _{B(t)}
\Delta \log \vert \vert J(\varphi) \vert \vert _\omega ^2e^{f_\varepsilon \circ \varphi}d\lambda\geq
C_\ep
\int_0^r{{dt}\over {t^{2n-1}}}\int _{B(t)}\vert \vert J(\varphi) \vert \vert _\omega ^{2/n}d\lambda
\leqno (19)$$

\noindent Maintenant, on rappelle la formule de Jensen suivante:

$$\int_0^r{{dt}\over {t^{2n-1}}}\int _{B(t)}
\Delta \rho d\lambda= \int_{\vert z \vert = r}\rho d\sigma- \rho (0)$$
valable pour toute fonction $\rho$ assez r\'eguli\`ere 
pour que les quantit\'es ci-dessus aient un sens.

\noindent Gr\^ace \`a cette identit\'e, le terme se trouvant de gauche de 
l'in\'egalit\'e (19) est \'egal \`a $\displaystyle \int_{\vert z \vert = r}
\log \vert \vert J(\varphi) \vert \vert ^2d\sigma + {\cal O}(1)$.
La concavit\'e de la fonction logarithme et les consid\'erations pr\'ec\'edentes
impliquent l'in\'egalit\'e suivante:

$$\log \int_{\vert z \vert = r}\vert \vert J(\varphi) \vert \vert_\omega ^{2/n}d\sigma
+ {\cal O}(1)\geq C 
\int_0^r{{dt}\over {t^{2n-1}}}\int _{B(t)}\vert \vert J(\varphi) \vert \vert_\omega ^{2/n}d\lambda
\leqno (20)$$

\noindent La preuve sera achev\'ee si l'on montre que l'in\'egalit\'e
(20) est impossible \`a satisfaire, pour des valeurs de $r$ assez grandes.

Pour chaque $r> 0$, soit $\displaystyle 
{\cal J}(r):= \int_0^r{{dt}\over {t^{2n-1}}}
\int _{B(t)}\vert \vert J(\varphi) \vert \vert _\omega ^{2/n}d\lambda$. La fonction ainsi d\'efinie est
croissante, et en d\'erivant succesivement on obtient:

$$r^{2n- 1}\int_{\vert z \vert = r}\vert \vert J(\varphi) \vert \vert ^{2/n}d\sigma=
\bigl (r^{2n- 1}{\cal J}^\prime\bigr)^\prime$$

\noindent pour toute valeur de $r$. On fait appel maintenant au lemme suivant 
(du type E. Borel), classique dans la th\'eorie de Nevanlinna (voir [7]).

\medskip
\noindent {\bf Lemme. }{\sl Soit $F:[0, \infty)\rightarrow \bR_+$ croissante, d\'erivable.
Alors pour tout $\delta> 0$, il existe un ensemble $E_\delta\subset \bR_+$
tel que $\displaystyle \int_{E_\delta}d\log t<\infty$ et tel que l'on ait:
l'in\'egalit\'e $rF^\prime(r)\leq F^{1+ \delta}(r)$,
si $r\in \bR_+\setminus E_\delta$.}
\medskip

\noindent Dans notre situation, on applique d'abord l'enonc\'e pr\'ec\'edent \`a la fonction
$F(r)= r^{2n-1}{\cal J}^\prime (r)$, avec $\delta= 1/2n-1$. Le lemme pr\'ec\'edent 
montre que 

$$r\bigl (r^{2n-1}{\cal J}^\prime (r)\bigr)^\prime\leq 
\bigl(r^{2n-1}{\cal J}^\prime (r)\bigr)^{1+ 1/2n-1}= r^{2n}
\bigl({\cal J}^\prime (r)\bigr)^{1+ 1/2n-1}.$$

\noindent Une nouvelle application du lemme, cette fois \`a la fonction
$F(r)= {\cal J}(r)$, montre que ${\cal J}^\prime (r)\leq {\cal J}^2$, en dehors
d'un ensemble de mesure logarithmique finie. En conclusion, 
on aura 

$$\bigl (r^{2n-1}{\cal J}^\prime (r)\bigr)^\prime\leq r^{2n-1}{\cal J}^4(r)\leqno (21)$$

\noindent pour tout $r\in \bR_+\setminus E$. 
On r\'e-ecrit maintenant l'in\'egalit\'e (15) comme suit:

$$4\log\int_0^r{{dt}\over {t^{2n-1}}}
\int _{B(t)}\vert \vert J(\varphi) \vert \vert _\omega ^{2/n}d\lambda+ {\cal O}(1)\leq 
C\int_0^r{{dt}\over {t^{2n-1}}}
\int _{B(t)}\vert \vert J(\varphi) \vert \vert _\omega  ^{2/n}d\lambda$$

\noindent si $r\in \bR_+\setminus E$. On choisit maintenant une suite 
$(r_k)\subset \bR_+\setminus E$ telle que $r_k\mapsto\infty$.  La derni\`ere relation donne alors une contradiction.

\vskip 30pt

\noindent {\twelvebf Bibliographie}

{
\parskip=4pt plus 1pt minus 1pt

\vskip 10pt

\item {[1]} Boucksom, S. {\sl On the volume of a line bundle} AG/0201031 (2002).
\smallskip 
\item {[2]} Boucksom, S. {\sl Divisorial Zariski decomposition } AG/0204336 (2002).
\smallskip 
\item {[3]} Boucksom, S., Demailly, J.-P., P\u aun, M., Peternell, T. 
{\sl The pseudo-effective cone of a compact K\"ahler manifold and varieties of 
negative Kodaira dimension}, AG/.
\smallskip 
\item{[4]} Campana, F. Connexit\'e rationnelle des vari\'et\'es de Fano. Ann. Sc. ENS. 25 (1992), 539-545.

\smallskip
\item {[5]} Demailly, J.-P. {\sl A numerical criterion for very ample line bundles}, 
J. Differential Geom 37 (1993) 323-374.
\smallskip 
\item {[6]} Demailly, J.-P. {\sl Regularization of closed positive 
currents of type (1,1) by the flow of a Chern connection}, Actes du Colloque en l'honneur de P. Dolbeault (Juin 1992), ŽditŽ par H. Skoda et J.M. TrŽpreau, Aspects of Mathematics, Vol. E 26, Vieweg, 
(1994) 105-126.
\smallskip 
\item {[7]} Griffiths, P. {\sl Entire holomorphic maps in one and several complex variables}
Annals of Math. Studies, Princeton Univ. Press., 1976. 
\smallskip 
\item {[8]} Kobayashi, S., Ochiai, T. {\sl Meromorphic mappings into compact complex spaces 
of general type}, Inv. Math. {\bf 31} (1975).
\smallskip 
\item {[9]} Kodaira, K. {\sl Holomophic mappings of polydiscs into compact complex manifolds}
J. Diff. Geom. {\bf 6} (1971).
\smallskip 

\smallskip

\item{[10]} Koll\`ar,J., Miyaoka, Y., Mori,S.{\sl Rationally connected Varieties}J.Alg.Geom. 1 (1992), 429-448.
\smallskip

\item {[11]} Lelong, P., Gruman, L. {\sl Entire Functions of Several Complex variables, } Sprin\-ger
1986.

\smallskip

\item{[12]} N. Sibony. Dynamique des applications rationnelles de $\b P^k$. in Panoramas et synth\`eses 8 (1999), 97-185. SMF. 

\smallskip 

\item {[13]} Tsuji, H. {\sl Pluricanonical systems of projective varieties of general type} preprint, 
AG/9909021.

\smallskip 

\item {[14]} Yau, S.-T. {\sl On the Ricci curvature of a 
complex K\"ahler manifold and the complex Monge--Amp\`ere equation},
Comm.\ Pure Appl.\ Math.\ {\bf 31} (1978).

}

\end